\documentclass[11pt]{article}
\usepackage[utf8]{inputenc}
\usepackage[T1]{fontenc}
\usepackage{lmodern}
\usepackage{geometry}
\usepackage{microtype}
\usepackage{hyperref}
\hypersetup{colorlinks=true,linkcolor=black,citecolor=black,urlcolor=blue}
\usepackage{setspace}
\usepackage{parskip} % space between paragraphs
\usepackage{csquotes}
\usepackage{enumitem} % for references list
\title{What makes a demonstration worthy of the name?}
\author{Paul Jorion\\\small ETHICS -- EA 7446, Universit\'e catholique de Lille\\\small 60 Bd Vauban, 59800 Lille, France}
\date{}
\begin{document}
\maketitle
\begin{abstract}
This paper revisits the foundations of mathematical proof through the lens of Aristotle’s threefold conception of truth: sensory evidence, axiomatic definition, and syllogistic deduction. I argue that modern mathematics has too often neglected the first of these—corroborated perception—in favour of axioms and theorems alone, thereby producing ``proofs'' that either disguise axioms as discoveries or repackage empirical observations in formal dress. Using Gödel’s incompleteness theorem as a case study, I show that Gödel’s reasoning covertly relies on an empirical notion of truth, smuggling physics into mathematics while claiming to remain purely formal. This reliance undermines the demonstration’s claim to be ``worthy of the name,'' since it oscillates between logical registers of differing validity—analytical, dialectical, and rhetorical. More broadly, the paper contends that mathematics often operates as a form of ``virtual physics,'' with intuition and empirical models quietly shaping its axioms and results. By clarifying the requirements for a genuine demonstration, I aim to restore the Aristotelian distinction between logic proper and looser modes of reasoning, and to highlight the risks of conflating mathematics with physics.
\end{abstract}

\section*{Three ways of ascertaining the truth}

One reads in an introductory book to the work of Kurt Gödel:

`` … Aristotelian logic rests on two pillars: a set of \textit{premises}, or \textit{assumptions}, which are taken to be true without proof, and a collection of \textit{rules of inference}, by which we transform one true statement into another […] As long as we can convince ourselves that the premises are true, then the conclusions [are] as solid and inescapable a fact as there ever will be. It is a certainty that follows from the semantic content of the premises and from the process of deduction laid down in our minds and formalised by Aristotle.

What Gödel discovered is that even if there exist true relations among pure numbers, the methods of deductive logic are just too weak for us to be able to prove all such facts. In other words, truth is simply bigger than proof '' (Casti \& DePauli 2000: 4-5).

Although widely accepted, what is asserted here in the second paragraph is based on a major misunderstanding as that truth is bigger than proof is no discovery of Gödel’s, as it was established by Aristotle twenty-five centuries ago and was the starting point of his groundbreaking reflection on truth and falsehood.

This fundamental misunderstanding, however, is far from trivial: it has undermined the pursuit of the mathematical venture in several ways. Indeed, mathematicians have attempted to prove statements which were in no need to be proved in the first place as their truth was established at little cost by other means, and in so doing they have at times introduced erroneous demonstrations aimed at proving… the undemonstrable, or hardly demonstrable (e.g. that there  exists in the Amazon a seven black-dotted yellow ladybird).

`` Truth is simply bigger than proof '', indeed, according to Aristotle, proof was but one of the three available sources of truth; since his days, none has been added to his brief list. Aristotle determined that there were three known ways for accessing the truth:

1° Truth known from perceptual evidence arising from our five senses, providing it can be corroborated so that facts are distinguished with certainty from illusions.

2° Truth derived from commonly accepted definitions, by which we assign to a single word or expression the meaning of a phrase: `` the fawn is the young of the deer and the doe ''; henceforth, whenever I am tempted to say `` the young of the deer and the doe '', I can say `` fawn '' instead. In mathematics, definitions are called `` axioms ''.

3° Truth reached as the conclusion of a well-formed syllogism, i.e. syntactically correct, by valid deductions from true premises. In mathematics, we call the equivalent of sequences of conclusions of valid syllogisms `` theorems ''.

Transposed into the language of mathematics, these three modes of accessing the truth become

1° Truth resulting from observations of the physical world expressed as models in mathematical language (mathematics as `` virtual physics '')

2° Truth established by axioms, which are definitions expressed in mathematical language

3° Truth obtained by the demonstration of a theorem, which is the conclusion of the proof provided.

Because in recent centuries mathematicians have ignored corroborated sensory evidence as the main and primary source of truth, they have attempted to prove by means of axioms and theorems solely, a considerable number of things that appear to be true, neglecting the `` obvious '' truths that come to us as trusted perceptions by our senses, i.e. as data captured from the world surrounding us as it stands. E.g. \textit{Goodstein}\textit{’}\textit{s theorem} which `` proves '' … an observation about natural numbers. Thus, instead of saying: `` Here is a curiosity that arises when we define a sequence of natural numbers in a certain way '', Goodstein thought it appropriate to demonstrate a theorem stating what we observe, taking axioms as the starting point for his demonstration. The fact that \textit{Goodstein}\textit{’}\textit{s theorem} can be demonstrated within one mathematical framework (second-order arithmetic), but not in a closely related one (Peano’s arithmetic), betrays the contrived nature of trying to mathematically prove a fact of observation.

In their desire to prove everything, mathematicians have produced erroneous proofs: clumsy and error-ridden attempts to prove what was in no need to be demonstrated and was often essentially unprovable due to the complexity of a demonstration candidate. Indeed, most of our surrounding empirical world requires in order to be proven an overwhelming number of elements requiring to be proven upstream. The mere fact of an observation as such summons the existence of an observer who would need to be proven in the first place.

Since sensory evidence has done the job satisfactorily in most instances, this wealth of unnecessary proofs were achieved through resorting to various tricks such as starting from false premises (paradoxes, conclusions of fallacious reasonings), or presenting as the conclusion of a demonstration worthy of the name some truth being in disguise either corroborated sensory evidence translated into mathematical language (\textit{virtual physics}) or an axiom. Such a `` theorem-in-no-need-to-be-demonstrated '' happens then to be true

1° either because it is in fact an axiom in disguise (i.e. a definition) where the demonstration is a \textit{petitio principii}: a circular reasoning, where what one claims to have discovered in fine was known from the start, the demonstration having been but an artful reformulation of the axiom

2° or because an observation derived from sensory evidence, an \textit{empirical fact}, is hidden in one of the axioms; in such a case, mathematicians have used mathematics to build a `` virtual physics '' and it is their intuition of the world around us that has guided them unbeknownst.

In one particular instance what is usually the `` stowaway '' presence of empirical facts within mathematical reasoning has been deliberate and proclaimed, I’m referring here to Hilbert’s \textit{6th proposition}:

\begin{quote}
\textbf{`` }\textit{Mathematical Treatment of the Axioms of Physics}. The investigations on the foundations of geometry suggest the problem: To treat in the same manner, by means of axioms, those physical sciences in which already today mathematics plays an important part; in the first rank are the theory of probabilities and mechanics. ''
\end{quote}

I reported in \textit{Comment la v}\textit{é}\textit{rit}\textit{é }\textit{et la r}\textit{é}\textit{alit}\textit{é }\textit{furent invent}\textit{é}\textit{es} (2009) that in the historical case of the invention of the infinitesimal calculus, the central issue was not of the protagonists at work, Newton’s and Leibniz’ intuition, whose workings would have remained opaque to their conscious self but of their quite explicit experimental validation all along of the mathematics whose methodology they were developing. Indeed, they defined the rules of this new method of calculation in such a way that their model of a planet’s orbit accurately reflected the observed motion, revealing the pre-eminence they implicitly gave to physics over mathematics in the development of their new tool (Jorion 2009: 332-346). The alternative approach would have been to demand mathematical rigour above all, which would have betrayed slight anomalies in the orbit of the planets, the pre-eminence being given in this case to mathematical precepts rather than to an accurate account of physical phenomena.

\textbf{\textit{What requirements for a demonstration to be }}\textbf{\textit{`` }}\textbf{\textit{worthy of the name}}\textbf{\textit{ ''}}

In my analysis of Gödel’s demonstration of his incompleteness theorem of arithmetic in \textit{Comment la v}\textit{é}\textit{rit}\textit{é }\textit{et la r}\textit{é}\textit{alit}\textit{é }\textit{furent invent}\textit{é}\textit{es} (2009: 285-326), I show in particular how at a critical juncture Gödel cheats by surreptitiously slipping into his demonstration that a proposition is true when its truth comes neither from being an axiom nor a demonstrated theorem, but from resulting from empirical evidence. In other words, at the moment when he finds himself mired in a paradox of his own construction, Gödel sneaks physics into mathematics in order to extricate himself from the pit he’s dug himself in.

Gödel’s proof of his incompleteness theorem of arithmetic fails to be `` worthy of the name '' for two types of reasons: for using arguments (\textit{transformation rules}) of too low a level of validity (\textit{demonstrability}) at different stages of his proof and for relying on an intuitive notion of truth whose precise concept varies according to his whim in different contexts.

Instead of sticking to arguments deemed by Aristotle to be of an `` analytical '' standard of logic, i.e. worthy of being used in a true demonstration with `` scientific '' status, Gödel went astray by resorting to arguments deemed more modestly to be of a `` dialectical '' standard of logic, i.e. acceptable in the courts or within the political debate, such as the \textit{reductio ad absurdum}, or even going so low as to use purely `` rhetorical '' arguments, i.e. acceptable only in everyday conversation. Of only a \textit{rhetorical} standard in Gödel’s demonstration is the recourse to a proposition asserting about itself that it is unprovable, having no other evidence to support this assertion than `` its own word '', the near equivalent in everyday conversation of `` It’s true because I’m telling you so ''. Attention is often drawn to the fact that the weakness here in the demonstration is due to the assertion being self-referential, but the bottom line is that there is nothing to support this pseudo-affirmation of a proposition about itself: a proposition cannot be `` taken at its word '' because there is no way a metamathematical argument can be held accountable and shamed for having lied.

When Gödel says of that proposition, as we’ve seen, that it is true even though it is impossible to prove this means out of necessity that it can only be true for one of the other two possible reasons why a proposition may be true: either that it is an axiom (i.e. a definition) or that it derives from corroborated sensory evidence. Since the context of the demonstration excludes the proposition in question from being a vast tautology: a hidden or unrecognised axiom, the only remaining option is that its truth derives from sensory evidence. But this would imply that the proposition has as its reference the empirical world around us, offering of it a true account, which would mean that Gödel is unwittingly bringing physics to bear on what he claims is mathematics only.

How is it possible that Gödel did not notice such a confusion in his argument between mathematics and physics? Because it is common knowledge that for him the two are coalesced: Gödel was a Platonist by ideological creed, convinced that the world is made up of numbers in its intimate build-up and that, as one of the consequences of this belief, everything we know in the world by sense evidence can also be proved mathematically. Bertrand Russell said of him: `` Gödel turned out to be an unadulterated Platonist, and apparently believed that an eternal ‘not’ was laid up in heaven, where virtuous logicians might hope to meet it hereafter '' (in Russell’s \textit{Autobiography}, quoted by Dawson 1988: 8).

I am arguing here, to the contrary, that to be `` worthy of the name '', a mathematical demonstration is not allowed to import surreptitiously elements of proof borrowed from physics: it is required to confine itself entirely to the formal domain of mathematics. It is easy to understand that were this distinction not upheld, the notion of an `` adequate mathematical model of physical reality '' would lose all meaning as stowaway empirical evidence might be lurking in the model at undisclosed locations, blurring the distinction between which elements of the physical world are properly modelled and which have straightforwardly been imported, i.e. simply depicted instead of accounted for.

In their book quoted earlier, Casti \& DePauli, unintentionally betray their unease in this respect when they write: `` What Gödel discovered is that even if there exist true relations among pure numbers, the methods of deductive logic are just too weak for us to be able to prove all such facts. In other words, truth is simply bigger than proof. This fact does not seem too astonishing when put into the context of everyday life. We are all aware of things that we ‘know’, but that we feel can never be logically deduced in a formal, Aristotelian fashion. In fact, the Oxford don and philosopher J. L. Austin, when first told of Gödel’s result, remarked: ‘Who would have ever thought otherwise?’ It seems likely that the average person on the street would say the same thing if someone announced that not everything that’s true can be known by following a process of logical deduction. But not so for mathematicians! '' (2000 : 4-5). Casti \& DePauli should have written more specifically: `` But not so for mathematicians who are devotees of the Platonic creed '' whose belief is that mathematics and physics are one and the same, or at least don’t care distinguishing one from the other. Hence their indifference to the quality of the proof mobilised at various steps in their demonstrations, and Gödel’s specific legerdemain in this respect with his `` I can’t care less '' attitude.

That the crux of the issue, the source of all his confusion, lies in Gödel’s Platonic ideological commitment, Casti and DePauli bring it out even more – albeit accidentally – in a passage where they attempt to clarify his approach with the help of an illustration featuring… chocolate cakes:

`` Truths = all conceivable cakes satisfying the chocolate cake test.

Proof = all recipes for actually making chocolate cakes with the Chocolate Cake Machine.

Now comes the Big Question: Is there a recipe for every conceivable chocolate cake? Or, equivalently, is every true statement provable? What we’re asking here is whether there are ‘honest-to-god’ chocolate cakes in the Platonic universe of cakes for which no recipe can ever be given '' (ibid. 19-20).

`` In the Platonic universe… '', one could not have stated it any better: Gödel’s demonstration only makes sense and represents a meaningful progress in understanding on the stage of a theatre where the Real would be entirely made out of numbers and where, as an entailment, two representations are fused: that of the empirical world surrounding us, which is accounted for by physics, and that of the world of mathematical entities, tools that we devised in order to offer a useful, if only stylised, vision of the world around us.

\textbf{\textit{Mathematicians}}\textbf{\textit{’ }}\textbf{\textit{intuition and }}\textbf{\textit{`` }}\textbf{\textit{virtual physics}}\textbf{\textit{ ''}}

Mathematics allows us to design models of physical reality. Some of our models are causal: if one thing happens, another or others will follow. We can also do simulations: we know what elements of different types are involved in a phenomenon of a collective nature and how they interact: let’s act it out, let numerical instances of objects interact for a while and see what happens.

It is often in no time that new mathematical objects and methods are mobilised for new advances in physics. This is no accident: they get immediately tested by various researchers for possible use. Einstein, for example, turned to tensors as soon as that mathematical tool became available.

Barrow reports: `` The development of non-Euclidean geometry as a branch of pure mathematics by Riemann in the nineteenth century, and the study of mathematical objects called tensors was a godsend to the development of twentieth-century physics. Tensors are defined by the fact that their constituant pieces change in a very particular fashion when their coordinate labels are altered in completely arbitrary ways. This esoteric mathematical machinery proved to be precisely what was required by Einstein in his formulation of the general theory of relativity. Non-Euclidean geometry described the distortion of space and time in the presence of mass-energy, while the behaviour of tensors ensured that any law of nature written in tensor language would automatically retain the same form no matter what the state of motion of the observer. Indeed, Einstein was rather fortunate in that his long-time friend, the pure mathematician Marcel Grossmann, was able to introduce him to these mathematical tools. Had they not already existed, Einstein could not have formulated the general theory of relativity '' (Barrow [1990] 1992: 189)

Not all of mathematics lends itself however with equal ease to the construction of physical models. Is there something then about a new mathematical tool that predisposes it to be summoned for modelling physical phenomena?

How are axioms defined? In \textit{Comment la v}\textit{é}\textit{rit}\textit{é }\textit{et la r}\textit{é}\textit{alit}\textit{é }\textit{furent invent}\textit{é}\textit{es} I argued that they are quite often not entirely abstract and seem to display an inclination towards the easy generation of physical models. I pointed out that this was then no deliberate strategy, the mathematician author remaining in all likelihood unaware of the fact. When Stephen Wolfram writes `` Maybe there’s something special about the particular axioms used in mathematics. And certainly if one thinks they’re the ones that uniquely describe science and the world, there might be a reason for that '' (Wolfram 2020: 559), this is precisely what I’m having in mind when I talk of mathematics as `` virtual physics ''. By taking as a starting point axioms that are already de facto `` virtual physics '', theorems are generated that remain all along `` virtual physics '' so to say by construction.

When I mentioned the attempts of mathematicians to demonstrate what is in no need to be demonstrated because it is in truth a straightforward fact of observation, something that our senses have allowed to establish, I already pointed out that in certain cases, the axiom makes it possible to surreptitiously import into the theorem an element of physics: `` because hidden in one of the axioms is an observation derived from sensory evidence, an \textit{empirical fact}, is hidden in one of the axioms; in such a case, mathematicians have used mathematics to build a `` virtual physics '' and it is their intuition of the world around us that has guided them unbeknownst ''.

Curiously, if one thinks of the fact that they are working on formal objects whose syntax and semantics alone should be relevant, mathematicians readily admit that they rely on their intuition, which they suggest is a very abstract thing that is very difficult to pin down more precisely. For example, Alan Turing, in his 1938 thesis, \textit{Systems of Logic Based on Ordinals}, remarked: ‘The activity of intuition consists in the production of spontaneous judgments which do not result from conscious trains of reasoning […] I shall not attempt to explain this idea of intuition more explicitly '' (in Copeland 2004: 192). In fact, the intuition of mathematicians owes its very shape to the world as it stands around us. And what mathematicians produce, therefore, since the origins of their discipline, and especially when they refer to their intuition as the source of their inspiration, is mathematics that constitutes, until better informed, a `` virtual physics ''.

This I set out systematically in \textit{Comment la v}\textit{é}\textit{rit}\textit{é }\textit{et la r}\textit{é}\textit{alit}\textit{é }\textit{furent invent}\textit{é}\textit{es}, taking as an example the birth of the infinitesimal calculus mentioned above, when it could be observed that the mathematics newly conceived by Newton and Leibniz were abused (the rigour asserted being sacrificed in the process) until they `` stuck '' to the celestial mechanics for which they were developed. The famous bishop \textit{cum} philosopher George Berkeley (1685-1753) denounced then the deception: a `` compensation of errors '' involving `` Ghosts of Departed quantities '' (Berkeley [1734] 1992: 199). He pointed out that `` in every other Science, Men prove their Conclusions [the ‘theorems’] by their Principles [the ‘axioms’], and not their Principles by their Conclusions. But if in yours you should allow your selves this unusual way of proceeding, the Consequence would be that you must take up with Induction and bid adieu to Demonstration. And if you submit to this, your Authority will no longer lead the way in Points of Reason and Science. I have no controversy about your Conclusions but only about your Logic and Method '' (ibid. 180). And Berkeley cruelly wondered `` whether such Mathematicians as cry out against [divine] Mysteries, have ever examined their own Principles? '' (ibid. 220).

You mathematicians, Berkeley argued, should take your starting point in mathematically meaningful axioms (`` Principles '') and see what theorem (`` Conclusions '') you can derive from them. But what you do instead is that you predefine the theorem (`` Conclusions ''), and then you make up the axioms (`` Principles '') that will allow you to prove it. In so doing, you do not proceed as scholars do, according to Aristotle, by \textit{deduction}, within the framework of the `` analytical '' standard of logic, but by \textit{induction}, which Aristotle has shown to be a mode of proof proper to the `` dialectictical '' standard of logic which is of choice for lawyers and politicians, who care only to not contradict themselves, but have no regard for the truth, unlike genuine scholars.

\section*{References}
\begin{itemize}[leftmargin=*,itemsep=0.2em]
  \item Barrow, John D., Theories of Everything. The Quest for Ultimate Explanation, London: Vintage 1990
  \item Berkeley, George, [1734] 1992 De Motu and The Analyst. A Modern Edition, with Introductions and Commentary, ed. and transl. from latin by Douglas M. Jesseph, Dordrecht (Holland): Kluwer Academic Publishers
  \item Casti, John L. \& Werner DePauli, Gödel. A Life of Logic, Cambridge (Mass.) Perseus: 2000
  \item Copeland, B. Jack, The Essential Turing, Oxford: Oxford University Press 2004
  \item Dawson, John W., `` Kurt Gödel in Sharper Focus '', in S. G. Shanker, Gödel’s Theorem in Focus, London: Croom Helm 1988
  \item Jorion, Paul, Comment la vérité et la réalité furent inventées, Bibliothèque des sciences humaines, Paris: Gallimard 2009
  \item Wolfram, Stephen, A Project to Find the Fundamental Theory of Physics, Wolfram Media 2020
\end{itemize}

\end{document}